\documentclass[leqno,oneside,english,12]{article}
\usepackage{amssymb,amsmath}

\author{Mariam Dhayni \thanks{Universit\'e d'Angers, Math\'ematiques,
49045 Angers ceded 01, France, e-mail:dhayni@math.univ-angers.fr}}
\title{Wilf's Conjecture for numerical semigroups
}
\date{\mbox{}}

\newtheorem{teorema}{Theorem}[section]

\newtheorem{proposicion}[teorema]{Proposition}
\newtheorem{lema}[teorema]{Lemma}

\newtheorem{corolario}[teorema]{Corollary}
\newtheorem{nota}[teorema]{Remark}
\newtheorem{exemple}[teorema]{Example}

\newtheorem{theorema}[teorema]{Theorem}
\newenvironment{demostracion}{\noindent Proof}

\numberwithin{equation}{section}
\newcommand{\N}{\mathbb N}

\begin{document}
\maketitle

\noindent{\bf Abstract:} Let $S\subseteq \N$ be a numerical semigroup with multiplicity $m$, embedding dimension $\nu$ and conductor $c=f+1=qm-\rho$ for some $q,\rho\in\N$ with $\rho<m$. Let Ap$(S,m)=\{w_0<w_1<\ldots<w_{m-1}\}$ be the Ap\'ery set of $S$. The aim of this paper is to prove Wilf's Conjecture in some special cases. First, we prove that if $w_{m-1}\geq w_1+w_{\alpha}$ and $(2+\frac{\alpha-3}{q})\nu\geq m$ for some $1<\alpha<m-1$, then $S$ satisfies Wilf's Conjecture. Then, we prove the conjecture in the following cases: $(2+\frac{1}{q})\nu\geq m$, $m-\nu\leq 5$ and $m=9$. Finally, the  conjecture is proved if $w_{m-1}\geq w_{\alpha-1}+w_{\alpha}$ and $(\frac{\alpha+3}{3})\nu\geq m$ for some $1<\alpha<m-1$.
 

\section{Introduction and notations }\label{sec 1 introduction}
Let $\N$ denote the set of natural numbers, including $0$. A \textit{numerical semigroup} $S$ is an additive submonoid of $(\N,+)$ of finite complement in $\N$, that is $0\in S$, if $a,b\in S$ then $a+b\in S$ and $\N\setminus S$ is a finite set. The elements of $\N\setminus S$ are called the \textit{gaps} of $S$. The largest gap is denoted by $f=f(S)=$max$(\N\setminus S)$ and is called the \textit{Frobenius number} of $S$. The smallest non zero element $m=m(S)=$min$(S^*)$ is called the \textit{multiplicity} of $S$ and $ n=\vert\{s\in S,s<f(S)\}\vert$ is also denoted by $n(S)$. Every numerical semigroup $S$ is finitely generated, i.e. is of the form
$$ S=<g_1,\ldots,g_{\nu}>=\N g_1+\ldots+\N g_{\nu}$$
for suitable unique coprime integers $g_1,\ldots,g_{\nu}.$ The number of generators of $S$ is denoted by $\nu=\nu(S)$ and is called the \textit{embedding dimension} of $S$. An integer $x\in\N\setminus S$ is called a \textit{ pseudo-Frobenius number } if $x+S^*\subseteq S.$ The \textit{ type } of the semigroup, denoted by $t(S)$ is the cardinality of set of pseudo-frobenius numbers. The \textit{ Ap\'ery set} of $S$ with respect to $a\in S$ is defined as $ \text{Ap}(S,a)=\{s\in S;s-a\notin S\}$.
\newline

\noindent The \textit{Diophantine Frobenius Problem}, also called the \textit{Coin Problem}, is to find the largest number $f$ that cannot be written in the form $\sum_{i=1}^n a_ix_i$; $x_i \in \N$ for given coprime positive integers $a_1,\cdots,a_n $. This problem is related to the theory of numerical semigroups in the following way: let $S$ be the numerical semigroup generated by $g_1,\ldots,g_{\nu}$, then $f$ is simply the largest integer not belonging to $S$. Hence the problem is to find a formula for $f$ in terms of the set of generators of $S$.
For $\nu = 2$, the formula of $f(<g_1, g_2>)$ is given by Sylvester (\cite{Sylvester1}), for $\nu \geq 3$ the problem is much harder. It has been proved in (\cite{Ramrez1}), that in general $f(S)$ is not algebraic in the set of generators of $S$. In (\cite{Wilf1}) 1978 H. S. Wilf proposed a lower bound for the number of generators of $S$ in terms 
of the Frobenius number as follows:
$ f(S)+1\leq \nu(S)n(S).$
\newline

\noindent Although the problem has been considered by several authors (cf. \cite{Barucci1}, \cite{Amors1}, \cite{Dobbs1}, \cite{Eliahou}, \cite{Froberg1}, \cite{Kaplan1}, \cite{Sammartano1} , \cite{Zhai1} ), only special cases have been solved and it remains wide open. In (\cite{Dobbs1}), D. Dobbs and G. Matthews proved Wilf's Conjecture for $\nu\leq 3$. In (\cite{Kaplan1}) N. Kaplan proved it for $f+1\leq 2m$ and in (\cite{Eliahou}) S. Eliahou extended Kaplan's work for $f+1\leq 3m$.   

\medskip

\noindent This work is a generalisation of the case covered by A. Sammartano (\cite{Sammartano1}), who showed that Wilf's Conjecture holds for $2\nu\geq m$ and $m\leq 8$, based on the idea of counting the elements of $S$ in some intervals of length $m$. We use different intervals in order to get an equivalent form of Wilf's Conjecture and then to prove it in more cases. We also cover the case where $2\nu\geq m$.

\noindent Here are few more details on the contents of this paper. Section 2 is devoted to give some notations that will enable us in the same section to give an equivalent form of Wilf's Conjecture. Section 3 is the heart of the paper. Let $\text{Ap}(S,m)=\lbrace 0=w_0<w_1<\cdots <w_{m-1}\rbrace$. First, we show that Wilf's conjecture holds for numerical semigroups that satisfy $w_{m-1}\geq w_{1}+w_{\alpha}$ and $(2+\frac{\alpha-3}{q})\nu\geq m$ for some $1<\alpha<m-1$ where $f+1=qm-\rho$ for some $q\in\N$, $0\leq \rho\leq m-1$. Then we prove  Wilf's Conjecture for numerical semigroups with $m-\nu\leq 4$ in order to cover the case where $2\nu \geq m$, prove by Sammartano in  (\cite{Sammartano1}). We also show that a numerical semigroup with $m-\nu=5$ verify Wilf's Conjecture in order to prove the conjecture for $m=9$. Finally, we show in this section, using the previous cases, that Wilf's conjecture holds for numerical semigroups with $(2+\frac{1}{q})\nu\geq m$. In section 4 we prove Wilf's Conjecture for numerical semigroups with $w_{m-1}\geq w_{\alpha-1}+w_{\alpha}$ and $(\frac{\alpha+3}{3})\nu\geq m$ for some $1<\alpha<m-1$.

 \medskip
 
 \noindent A good reference on numerical semigroups is \cite{Rosale1}.

\section{Preliminaries}\label{sec 2 Notations and Equivalent form OfWilf's Conjecture}
Let the notations be as in the introduction. For the sake of clarity we shall use the notations $\nu,f,n,...$ for $\nu(S),f(S),n(S)...$.
In this section we will introduce some notations and family of numbers that will enable us to give an equivalent form of Wilf's conjecture.

\noindent {\bf Definition}{.} Let $S$ be a numerical semigroup and let $c=C(S)=f+1$ be the conductor of $S$. Denote by 
$$q=\lceil \frac{c}{m}\rceil,  $$
where $\lceil x\rceil$ denote the smallest integer greater than or equal to $x$.
Thus, $qm\geq c$ and $c=qm-\rho $ with $0\leq \rho <m$. 

\noindent Given a non negative integer $k$, we define the $k$th interval of length $m,$
$$I_k=[km-\rho,(k+1)m-\rho[=\{km-\rho,km-\rho+1,\ldots,(k+1)m-\rho-1\}.$$
 We denote by
$$n_k=\vert\{s\in S \cap I_k\}\vert.$$ 
For $j\in\{1,\ldots, m-1\}$, we define $\eta_j$ to be the number of intervals $I_k$ with  $n_k=j$.
$$ \eta_j=\vert\{k\in \N ; \vert I_k\cap S\vert=j\}\vert.$$
\begin{proposicion}\label{prop1}{\rm Under the previous notations, we have:
\begin{enumerate}
\item[$i)$] $1 \leq n_k \leq m-1$ for all $k=0,\ldots,q-1.$
\item[$ii)$] $n_k= m $ for all $k\geq q.$
\item[$iii)$] $n=n(S)= \sum_{k=0}^{q-1}n_k.$
\item[$iv)$]$\sum_{j=1}^{m-1}\eta_j=q.$
\item[$v)$]$\sum_{j=1}^{m-1}j\eta_j=\sum_{k=0}^{q-1} n_k=n$.
\end{enumerate}}
\end{proposicion}
\begin{demostracion}{.} $i),ii),iii)$ are obvious. Now using $i)$, $ii)$ we will prove $iv)$ and $v)$. 
\begin{enumerate}
\item[$iv)$]$\sum_{j=1}^{m-1}\eta_j= \sum_{j=1}^{m-1}\vert\{k \in \N; \vert I_k\cap S\vert=j\}\vert = \sum_{j=1}^{m-1}\vert\{k \in \N; n_k=j;k=0,\ldots,q-1\}\vert=q.$
\item[$v)$]$\sum_{j=1}^{m-1}j\eta_j = \sum_{j=1}^{m-1}j\vert\{k \in \N; \vert I_k\cap S\vert=j\}\vert = \sum_{j=1}^{m-1}j\vert\{k \in \N; n_k=j;k=0,\ldots,q-1\}\vert =\sum_{k=0}^{q-1} n_k=n.   \quad \square
$
\end{enumerate}
\end{demostracion}

\medskip

\noindent \textbf{Remark:} {\rm We shall use the notation $\lfloor x\rfloor$ for the largest integer smaller than or equal to $x$.}

\bigskip

\noindent Next we will express $\eta_j$ in terms of th Ap\'ery set.

\begin{proposicion}\label{prop3}\rm Let Ap$(S,m)=\{w_0=0<w_1<w_2<\ldots <w_{m-1}\}$. Under the previous notations, we have for all $j\in\{1,\ldots,m-1\}$
\begin{equation*}
\eta_j=\lfloor\frac{w_j+\rho}{m}\rfloor-\lfloor \frac{w_{j-1}+\rho}{m}\rfloor.
\end{equation*}
\end{proposicion}
\begin{demostracion}{.} Fix $0\leq k\leq q-1$,  and let $ j\in \{1,\ldots,m-1\}$. We will show that the interval $I_k$ contains exactly $j$ elements of $S$ if and only if $w_{j-1}<(k+1)m-\rho\leq w_j$.  
Recall to this end that  for all $s\in S$, there exist $0\leq i\leq m-1,a\in \N$ such that $s=w_i+am$.

\medskip

\noindent Suppose that $I_k$ contains $j$ elements. Suppose, by contradiction, that $w_{j-1} \geq (k+1)m-\rho$. We have $w_{m-1}>\ldots>w_{j-1} \geq (k+1)m-\rho $, thus $w_{m-1},\ldots,w_{j-1}\in \cup_{t=k+1}^qI_t$. Hence, $I_k$ contains at most $j-1$ elements of $S$ (namely $w_0+km=km,w_1+k_1m,w_2+k_2m,\ldots,w_{j-2}+k_{j-2}m$ for some $k_1,\ldots,k_{j-2}\in \{0,\ldots,k-1\}$). This contradicts the fact that $I_k$ contains exactly $j$ elements of $S$.
\\
\noindent If $w_j<(k+1)m-\rho$, then $w_0<\ldots<w_j<(k+1)m-\rho$, thus $w_0,\ldots, w_j\in \cup_{t=0}^kI_t$. Hence, $I_k$ contains at least $j+1$ elements of $S$ which are : $w_0+km=km,w_1+k_1m,w_2+k_2m,\ldots,w_{j}+k_{j}m$ for some $k_1,\ldots,k_{j}\in \{0,\ldots,k-1\}$, which contradicts the fact that $I_k$ contains exactly $j$ elements of $S$. 

\medskip

\noindent Conversely, $w_{j-1}<(k+1)m-\rho$ implies that $w_0<\ldots< w_{j-1}<(k+1)m-\rho$, so  $w_0,\ldots, w_{j-1}\in \cup_{t=0}^kI_t$. Thus $I_k$ contains at least $j$ elements which are : $w_0+km=km,w_1+k_1m,w_2+k_2m,\ldots,w_{j-1}+k_{j-1}m$ for some $k_1,\ldots,k_{j-1}\in \{0,\ldots,k-1\}$.
\\
\noindent On the other hand  $w_{j} \geq (k+1)m-\rho$ implies that $w_{m-1}>\ldots>w_{j} \geq (k+1)m-\rho $, so $w_{m-1},\ldots,w_{j}\in \cup_{t=k+1}^qI_t$. Thus $I_k$ contains at most $j$ elements which are : $w_0+km=km,w_1+k_1m,w_2+k_2m,\ldots,w_{j-1}+k_{j-1}m$ for some $k_1,\ldots,k_{j-1}\in \{0,\ldots,k-1\}$. Hence, if  $w_{j-1}<(k+1)m-\rho\leq w_j$, then $I_k$ contains exactly $j$ elements of $S$ and this proves our assertion.\\
 Consequently,
\medskip
\begin{equation*}
\begin{array}{rcl}
\eta_j & =& \left| \{k\in \N \mbox{ such that } \vert I_k \cap S\vert =j\}\right|\\ \\
 & =& \vert\{ k\in \N \mbox{ such that } w_{j-1}<(k+1)m-\rho\leq w_j\}\vert\\ \\ 
  & =& \vert\{ k\in \N \mbox{ such that } \frac{w_{j-1}+\rho}{m}<(k+1)\leq \frac{w_j+\rho}{m}\}\vert\\ \\
& =& \vert\{ k\in \N \mbox{ such that } \frac{w_{j-1}+\rho}{m}-1<k\leq \frac{w_j+\rho}{m}-1\}\vert\\ \\
& =& \vert\{ k\in \N \mbox{ such that } \lfloor\frac{w_{j-1}+\rho}{m}\rfloor\leq k\leq \lfloor\frac{w_j+\rho}{m}\rfloor-1\}\vert\\ \\
 & = & \lfloor\frac{w_j+\rho}{m}\rfloor-\lfloor\frac{w_{j-1}+\rho}{m}\rfloor. \quad \square 
\end{array}
\end{equation*} 
\end{demostracion} 
\noindent Proposition \ref{prop2} gives an equivalent form of Wilf's Conjecture using  Propositions \ref{prop1} and \ref{prop3}.
\begin{proposicion}\label{prop2}{\rm Let $S$ be a numerical semigroup with multiplicity $m$, embedding dimension $\nu$ and conductor $f+1=qm-\rho$ for some $q\in\N$ and $0\leq\rho\leq m-1$. Let $w_0=0<w_1<w_2<\ldots <w_{m-1}$ be the elements of Ap$(S,m)$. Then $S$ satisfies Wilf's Conjecture  \textit{if and only if} }
\begin{equation*}
\sum_{j=1}^{m-1}(\lfloor\frac{w_j+\rho}{m}\rfloor-\lfloor\frac{w_{j-1}+\rho}{m}\rfloor)(j\nu-m)+\rho\geq 0.
\end{equation*}
\end{proposicion}
\begin{demostracion}{.} By Proposition \ref{prop1}, we have
\begin{equation*}
  f+1\leq n\nu \Leftrightarrow  qm-\rho  \leq  \nu \sum_{k=0}^{q-1}n_k  
 \Leftrightarrow \sum_{k=0}^{q-1}m-\rho  \leq \sum_{k=0}^{q-1}n_k\nu
 \Leftrightarrow \sum_{k=0}^{q-1}(n_k\nu-m)+\rho  \geq  0\Leftrightarrow \sum_{j=1}^{m-1}\eta_j(j\nu-m)+\rho \geq 0.
\end{equation*}
Using Proposition \ref{prop3}, we have
\begin{equation*}
\sum_{j=1}^{m-1}\eta_j(j\nu-m)+\rho \geq 0 \Leftrightarrow \sum_{j=1}^{m-1}(\lfloor\frac{w_j+\rho}{m}\rfloor-\lfloor\frac{w_{j-1}+\rho}{m}\rfloor)(j\nu-m)+\rho\geq 0.
\end{equation*}
Thus the proof is complete. $ \quad \square $

\end{demostracion}
\begin{nota}\label{nota1}
{\rm Let Ap($S,m)=\{w_0=0<w_1<\ldots<w_{m-1}\}$. The following technical results will be used through the paper:
\begin{enumerate}
\item $\lfloor\frac{w_0+\rho}{m}\rfloor =0$ ($w_0=0$ and $0\leq\rho<m$).
\item For all $1\leq i\leq m-1$, $\lfloor\frac{w_i+\rho}{m}\rfloor\geq 1$ ($w_i>m$).
\item  For all $1\leq i\leq m-1$, either $\lfloor\frac{w_i+\rho}{m}\rfloor= \lfloor\frac{w_i}{m}\rfloor$ or $\lfloor\frac{w_i+\rho}{m}\rfloor=\lfloor\frac{w_i}{m}\rfloor+1$. In the second case $\lfloor\frac{w_i+\rho}{m}\rfloor\geq 2$ and $\rho\geq 1$.
\item If $w_i<w_j$ for some $0\leq i<j\leq m-1$, then $\lfloor\frac{w_i+\rho}{m}\rfloor\leq\lfloor\frac{w_j+\rho}{m}\rfloor.$
\item $\lfloor\frac{w_{m-1}+\rho}{m}\rfloor=\lfloor\frac{qm-\rho+m-1+\rho}{m}\rfloor=q$.
\end{enumerate}
}

\end{nota}

\section{Main Results}\label{sec 4 N.S }
In this section, we show that Wilf's Conjecture holds for numerical semigroups in the following cases: 
\begin{enumerate}
\item $w_{m-1}\geq w_1+w_{\alpha}$ and $(2+\frac{\alpha-3}{q})\nu\geq m$ for some $1<\alpha<m-1$.
\item $m-\nu\leq 5$. (Note that the case $m-\nu\leq 4$ results from the fact that Wilf's Conjecture holds for $2\nu \geq m$ (\cite{Sammartano1}), however we shall give the proof for $m-\nu\leq 3$ in order to cover this result through our techniques).
\end{enumerate}
\noindent We then deduce the conjecture for $m=9$ and for $(2+\frac{1}{q})\nu\geq m$.\\

\medskip

\noindent The following technical Lemma will be used through the paper:

\begin{lema}\label{x}
{\rm Let Ap($S,m)=\{w_0=0<w_1<\ldots<w_{m-1}\}$. Suppose that $w_i\geq w_j+w_k$, then $\lfloor\frac{w_i+\rho}{m}\rfloor\geq \lfloor\frac{w_j+\rho}{m}\rfloor+\lfloor\frac{w_k+\rho}{m}\rfloor-1$. If furthermore, $\lfloor\frac{w_i+\rho}{m}\rfloor-\lfloor\frac{w_j+\rho}{m}\rfloor-\lfloor\frac{w_k+\rho}{m}\rfloor=-1$, then $\lfloor\frac{w_j+\rho}{m}\rfloor=\lfloor\frac{w_j}{m}\rfloor+1$, $\lfloor\frac{w_k+\rho}{m}\rfloor=\lfloor\frac{w_k}{m}\rfloor+1$ and $\rho\geq 1$. In particular, $\lfloor\frac{w_j+\rho}{m}\rfloor\geq 2$, $\lfloor\frac{w_k+\rho}{m}\rfloor\geq 2$ and $\rho\geq 1$.}
\end{lema}

\begin{demostracion}{.} $w_i\geq w_j+w_k$ implies that $w_i+\rho\geq w_j+w_k+\rho$, hence 
$\frac{w_i+\rho}{m}\geq\frac{ w_j+w_k+\rho}{m}.$ Consequently, $\lfloor\frac{w_i+\rho}{m}\rfloor\geq\lfloor\frac{ w_j+w_k+\rho}{m}\rfloor.$ Therefore, $\lfloor\frac{w_i+\rho}{m}\rfloor\geq\lfloor\frac{ w_j+\rho}{m}\rfloor+\lfloor\frac{ w_k}{m}\rfloor.$ Hence, by Remark \ref{nota1} (3), $\lfloor\frac{w_i+\rho}{m}\rfloor\geq \lfloor\frac{w_j+\rho}{m}\rfloor+\lfloor\frac{w_k+\rho}{m}\rfloor-1$.\\
\noindent Suppose that  $w_i\geq w_j+w_k$ and that $\lfloor\frac{w_i+\rho}{m}\rfloor-\lfloor\frac{w_j+\rho}{m}\rfloor-\lfloor\frac{w_k+\rho}{m}\rfloor=-1$. Suppose by the way of contradiction that either $\lfloor\frac{w_j+\rho}{m}\rfloor\neq\lfloor\frac{w_j}{m}\rfloor+1$  or $\lfloor\frac{w_k+\rho}{m}\rfloor\neq\lfloor\frac{w_k}{m}\rfloor+1$ or $\rho<1$. Then, by Remark \ref{nota1} (3) and the fact that $\rho\geq 0$, we have either $\lfloor\frac{w_j+\rho}{m}\rfloor=\lfloor\frac{w_j}{m}\rfloor$ or $\lfloor\frac{w_k+\rho}{m}\rfloor=\lfloor\frac{w_k}{m}\rfloor$ or $\rho=0$. Since $w_i\geq w_j+w_k$, then $\lfloor\frac{w_i+\rho}{m}\rfloor\geq\lfloor\frac{ w_j+w_k+\rho}{m}\rfloor.$ In this case $\lfloor\frac{w_i+\rho}{m}\rfloor\geq\lfloor\frac{ w_j+\rho}{m}\rfloor+\lfloor\frac{ w_k+\rho}{m}\rfloor,$ which is impossible. Hence, $\lfloor\frac{w_j+\rho}{m}\rfloor=\lfloor\frac{w_j}{m}\rfloor+1$, $\lfloor\frac{w_k+\rho}{m}\rfloor=\lfloor\frac{w_k}{m}\rfloor+1$ and $\rho\geq 1$. Therefore, by Remark \ref{nota1} (2), $\lfloor\frac{w_j+\rho}{m}\rfloor=\lfloor\frac{w_j}{m}\rfloor+1\geq2$, $\lfloor\frac{w_k+\rho}{m}\rfloor=\lfloor\frac{w_k}{m}\rfloor+1\geq2$ and $\rho\geq 1$. $\quad \square$

\end{demostracion}

\medskip
\noindent Next we will show that Wilf's Conjecture holds for numerical semigroups with $w_{m-1}\geq w_1+w_{\alpha}$ and $(2+\frac{\alpha-3}{q})\nu\geq m$.
\begin{theorema}\label{thm2} {\rm Let $S$ be a numerical semigroup with multiplicity $m$, embedding dimension $\nu$ and conductor $f+1=qm-\rho$ for some $q\in \N$, $0\leq \rho\leq m-1$. Let $w_0=0<w_1<w_2<\ldots<w_{m-1}$ be the elements of Ap($S,m)$. Suppose that $w_{m-1}\geq w_1+w_{\alpha}$ for some $1<\alpha<m-1$. If $(2+\frac{\alpha-3}{q})\nu\geq m$, then $S$ satisfies Wilf's Conjecture.
}
\end{theorema}
\begin{demostracion}{.} We are going to show that $S$ satisfies Wilf's Conjecture by means of Proposition \ref{prop2}. We have,

\medskip

\begin{equation}\label{d1}
\begin{array}{lll}
\displaystyle{\sum_{j={1}}^{\alpha}(\lfloor\frac{w_j+\rho}{m}\rfloor-\lfloor\frac{w_{j-1}+\rho}{m}\rfloor)}(j\nu-m)&=\displaystyle{\sum_{j={1}}^{\alpha}\lfloor\frac{w_j+\rho}{m}\rfloor(j\nu-m)-\sum_{j=1}^{\alpha}\lfloor\frac{w_{j-1}+\rho}{m}\rfloor}(j\nu-m)\\&=\displaystyle{\sum_{j={1}}^{\alpha}\lfloor\frac{w_j+\rho}{m}\rfloor(j\nu-m)-\sum_{j=0}^{\alpha-1}\lfloor\frac{w_j+\rho}{m}\rfloor}((j+1)\nu-m)\\&= \displaystyle{\lfloor\frac{w_{\alpha}+\rho}{m}\rfloor(\alpha\nu-m)-\lfloor\frac{w_{0}+\rho}{m}\rfloor(\nu-m)-\sum_{j=1}^{\alpha-1}\lfloor\frac{w_{j}+\rho}{m}\rfloor\nu}\\ &= \displaystyle{\lfloor\frac{w_{\alpha}+\rho}{m}\rfloor(\alpha\nu-m)-\lfloor\frac{w_{1}+\rho}{m}\rfloor\nu-\sum_{j=2}^{\alpha-1}\lfloor\frac{w_{j}+\rho}{m}\rfloor\nu}
\\& \geq \displaystyle{\lfloor\frac{w_{\alpha}+\rho}{m}\rfloor(\alpha\nu-m)-\lfloor\frac{w_{1}+\rho}{m}\rfloor\nu-\sum_{j=2}^{\alpha-1}\lfloor\frac{w_{\alpha}+\rho}{m}\rfloor\nu}
\\& = \displaystyle{\lfloor\frac{w_{\alpha}+\rho}{m}\rfloor(\alpha\nu-m)-\lfloor\frac{w_{1}+\rho}{m}\rfloor\nu-\lfloor\frac{w_{\alpha}+\rho}{m}\rfloor(\alpha-2)\nu}\\ \\ &= \displaystyle{-\lfloor\frac{w_{1}+\rho}{m}\rfloor\nu+\lfloor\frac{w_{\alpha}+\rho}{m}\rfloor(2\nu-m).}
\end{array}
\end{equation}

\begin{equation}\label{4.1 2}
\begin{array}{llll}
\displaystyle{\sum_{j={\alpha+1}}^{m-1}(\lfloor\frac{w_j+\rho}{m}\rfloor-\lfloor\frac{w_{j-1}+\rho}{m}\rfloor)(j\nu-m)}&\geq \displaystyle{\sum_{j={\alpha+1}}^{m-1}(\lfloor\frac{w_j+\rho}{m}\rfloor-\lfloor\frac{w_{j-1}+\rho}{m}\rfloor)((\alpha+1)\nu-m)}\\&=
\displaystyle{((\alpha+1)\nu-m)\sum_{j={\alpha+1}}^{m-1}(\lfloor\frac{w_j+\rho}{m}\rfloor-\lfloor\frac{w_{j-1}+\rho}{m}\rfloor)}\\

&=\displaystyle{((\alpha+1)\nu-m)(\sum_{j={\alpha+1}}^{m-1}\lfloor\frac{w_j+\rho}{m}\rfloor-\sum_{j=\alpha+1}^{m-1}\lfloor\frac{w_{j-1}+\rho}{m}\rfloor)}\\&=((\alpha+1)\nu-m) \displaystyle{(\sum_{j={\alpha+1}}^{m-1}\lfloor\frac{w_j+\rho}{m}\rfloor-\sum_{j=\alpha}^{m-2}\lfloor\frac{w_j+\rho}{m}\rfloor)}

\\&=\displaystyle{(\lfloor\frac{w_{m-1}+\rho}{m}\rfloor-\lfloor\frac{w_{\alpha}+\rho}{m}\rfloor)((\alpha+1)\nu-m).}
\end{array}
\end{equation}

\noindent Since $w_{m-1}\geq w_1+w_{\alpha}$, by Lemma \ref{x}, it follows that $\lfloor\frac{w_{m-1}+\rho}{m}\rfloor\geq\lfloor\frac{w_1+\rho}{m}\rfloor+\lfloor\frac{w_{\alpha}+\rho}{m}\rfloor-1.$ Let $x=\lfloor\frac{w_{m-1}+\rho}{m}\rfloor-\lfloor\frac{w_1+\rho}{m}\rfloor-\lfloor\frac{w_{\alpha}+\rho}{m}\rfloor.$ Then, 
$\lfloor\frac{w_1+\rho}{m}\rfloor+\lfloor\frac{w_{\alpha}+\rho}{m}\rfloor=q-x \text{ and } x\geq -1$. Now using \eqref{d1} and \eqref{4.1 2}, we have\\
\begin{equation*}
\begin{array}{lll}
\displaystyle{\sum_{j=1}^{m-1}(\lfloor\frac{w_j+\rho}{m}\rfloor-\lfloor\frac{w_{j-1}+\rho}{m}\rfloor)(j\nu-m)+\rho}
&\displaystyle{\geq -\lfloor\frac{w_1+\rho}{m}\rfloor\nu+ \lfloor\frac{w_{\alpha}+\rho}{m}\rfloor(2\nu-m)}\\ &
\displaystyle{+(\lfloor\frac{w_{m-1}+\rho}{m}\rfloor-\lfloor\frac{w_{\alpha}+\rho}{m}\rfloor)((\alpha+1)\nu-m)+\rho}
\\ \\&\displaystyle{=\lfloor\frac{w_1+\rho}{m}\rfloor(-\nu+((\alpha+1)\nu-m ) -( (\alpha+1)\nu-m))    + \lfloor\frac{w_{\alpha}+\rho}{m}\rfloor(2\nu-m)}\\ &
\displaystyle{+(\lfloor\frac{w_{m-1}+\rho}{m}\rfloor-\lfloor\frac{w_{\alpha}+\rho}{m}\rfloor)((\alpha+1)\nu-m)+\rho}
\\ \\
&\displaystyle{=\lfloor\frac{w_1+\rho}{m}\rfloor(\alpha\nu-m)+\lfloor\frac{w_\alpha+\rho}{m}\rfloor(2\nu-m)}\\ &\displaystyle{+(\lfloor\frac{w_{m-1}+\rho}{m}\rfloor-\lfloor\frac{w_{\alpha}+\rho}{m}\rfloor-\lfloor\frac{w_1+\rho}{m}\rfloor)((\alpha+1)\nu-m)+\rho}\\ \\
&\displaystyle{=(\lfloor\frac{w_1+\rho}{m}\rfloor+\lfloor\frac{w_\alpha+\rho}{m}\rfloor)(2\nu-m)}+\lfloor\frac{w_1+\rho}{m}\rfloor(\alpha-2)\nu\\ &\displaystyle{+(\lfloor\frac{w_{m-1}+\rho}{m}\rfloor-\lfloor\frac{w_{\alpha}+\rho}{m}\rfloor-\lfloor\frac{w_1+\rho}{m}\rfloor)((\alpha+1)\nu-m)+\rho}\\ \\
&\displaystyle{=(q-x)(2\nu-m)}+\lfloor\frac{w_1+\rho}{m}\rfloor(\alpha-2)\nu\displaystyle{+x((\alpha+1)\nu-m)+\rho}\\ \\
&\displaystyle{\geq(q-x)(2\nu-m)}+(\alpha-2)\nu\displaystyle{+x((\alpha+1)\nu-m)+\rho}\\ \\

&\displaystyle{=\nu(2q-2x+\alpha-2+x\alpha+x)-qm+\rho}\\ \\

&\displaystyle{=\nu(2q+(\alpha-2)(x+1)+x)-qm+\rho}\\ \\

&\displaystyle{\geq\nu(2q+\alpha-3)-qm+\rho} \hspace{4cm}(x\geq-1)\\ \\
&\displaystyle{=q(\nu(2+\frac{\alpha-3}{q})-m)+\rho}\\ \\
&\geq 0.
\end{array}
\end{equation*}

\noindent Using  Proposition \ref{prop2}, we get that $S$ satisfies Wilf's Conjecture. $\quad \square$
\end{demostracion}

\begin{exemple}\label{ex1} {\rm Consider the following numerical semigroup $S=<19,21,23,25,27,28>$. Note that $3\nu<m$. We have $w_1=21$, $w_{14}=56$ and $w_{m-1}=83$ i.e. $w_{m-1}\geq w_1+w_{14}$. In addition, $(2+\frac{\alpha-3}{q})\nu=(2+\frac{14                                                                                                                                                                                                                                                                                                                                                                                                                                                                                                                                                                                                                                                                                                                                                                                                                                                           -3}{4})6\geq19=m$. Thus the conditions of Theorem \ref{thm2} are valid, so $S$ satisfies Wilf's Conjecture.} 
\end{exemple}

\medskip
\noindent In the following we shall deduce some cases where Wilf's Conjecture holds. We start with the following technical Lemma.
\begin{lema}\label{lem1} {\rm Let $S$ be a numerical semigroup with multiplicity $m$ and embedding dimension $\nu$. Let $w_0=0<w_1<w_2<\ldots<w_{m-1}$ be the elements of Ap$(S,m)$. If $m-\nu> \displaystyle{(^\alpha_2)=\frac{\alpha(\alpha-1)}{2} }$ for some $\alpha\in\N^{*}$, then $w_{m-1}\geq w_1+w_{\alpha}$.
}
\end{lema}
\begin{demostracion}{.} Recall that an element $x$ of the Ap\'ery set of $S$ belongs to min(Ap($S,m))$ if and only if $x\neq w_i+w_j$ for all $w_i,w_j\in$Ap($S,m)\setminus\{0\}$, in particular $m-\nu=\vert$Ap$(S,m)\setminus$min(Ap($S,m))\vert$. Suppose by the way of contradiction that $w_{m-1}<w_1+w_{\alpha}$, and let  $w\in\text{Ap}(S,m)\setminus\min(\text{Ap}(S,m))$. Then $w\leq w_{m-1}$ and $w=w_i+w_j$ for some $w_i,w_j\in$ Ap$(S,m)\setminus\{0\}$. Hence,  $w\leq w_{m-1}< w_1+w_{\alpha}$. Thus the only possible values for $w$ are $\{w_i+w_j;1\leq i\leq j\leq \alpha-1\}$. Therefore, $m-\nu\leq\displaystyle{(^\alpha_2)=\frac{\alpha(\alpha-1)}{2} } $, which is impossible. Hence, $w_{m-1}\geq w_1+w_{\alpha}$.$\quad \square$
\end{demostracion}

\medskip

\noindent Next we will deduce Wilf's Conjecture for numerical Semigroups with $m-\nu> \frac{\alpha(\alpha-1)}{2}$ and $(2+\frac{\alpha-3}{q})\nu\geq m$. It will be used later to show that the conjecture holds for those with $(2+\frac{1}{q})\nu\geq m$, and inorder also to cover the result in \cite{Sammartano1} saying that the conjecture is true for $2\nu\geq m$.
\begin {corolario}\label{cor1}{\rm Let $S$ be a numerical semigroup with multiplicity $m$, embedding dimension $\nu$ and conductor $f+1=qm-\rho$ for some $q\in\N$, $0\leq \rho\leq m-1$. Suppose that  $m-\nu> (^\alpha_2)=\frac{\alpha(\alpha-1)}{2}$ for some $ 1<\alpha<m-1$. If $(2+\frac{\alpha-3}{q})\nu\geq m$, then $S$ satisfies Wilf's Conjecture.
}
\end{corolario}
\begin{demostracion}{.} It follows from Lemma \ref{lem1} that if $m-\nu>\frac{\alpha(\alpha-1)}{2}$, then $w_{m-1}\geq w_1+w_{\alpha}$. Now use Theorem \ref{thm2}. $\quad \square$

\end{demostracion}
\medskip
\noindent As a direct consequence of Theorem \ref{thm2}, we get the following Corollary.
\begin{corolario}\label{cor3}
{\rm Let $S$ be a numerical semigroup with a given multiplicity $m$ and conductor $f+1=qm-\rho$ for some $q\in\N$, $0\leq \rho\leq m-1$. Let $w_0=0<w_1<\ldots<w_{m-1}$ be the elements of Ap$(S,m)$. If $w_{m-1}\geq w_1+w_{\alpha}$ for some $ 1<\alpha<m-1$ and $m\leq 8+4(\frac{\alpha-3}{q})$ then $S$ satisfies Wilf's Conjecture.}
\end{corolario}
\begin{demostracion}{.}
 By Theorem \ref{thm2}, we may assume that $(2+ \frac{\alpha -3}{q})\nu<m$. Therefore, $\nu<\frac{qm}{2q+\alpha-3}\leq \frac{8q+\alpha-12}{2q+\alpha-3}.$ Hence $\nu<4$, consequently $S$ satisfies Wilf's Conjecture (\cite{Dobbs1}). $\quad \square$
\end{demostracion}

\medskip
\noindent In the following Lemma, we will show that Wilf's Conjecture holds for numerical semigroups with $m-\nu\leq 3$. This will enable us later to  prove the conjecture for numerical semigroups with $(2+\frac{1}{q})\nu\geq m$ and cover the result in \cite{Sammartano1} saying that the conjecture is true for $2\nu\geq m$.

\begin{lema}\label{thm1} \rm{Let $S$ be a numerical Semigroup with multiplicity $m$ and embedding dimension $\nu$. If $m-\nu\leq 3$, then $S$ satisfies Wilf's Conjecture.}
\end{lema}
\begin{demostracion}{.} We my assume that $\nu\geq 4$ ( $\nu\leq 3$ is solved \cite{Dobbs1}). We are going to show that $S$ satisfies Wilf's Conjecture by means of Proposition \ref{prop2}.
\begin{enumerate}
\item[\textbf{\textit{i})}] If $m-\nu=1$, then we may assume that $m=\nu+1\geq 5$.
By taking $\alpha=1$ in \eqref{4.1 2}, we get \newline
$
\begin{array}{lll}
\displaystyle{\sum_{j={2}}^{m-1}(\lfloor\frac{w_j+\rho}{m}\rfloor-\lfloor\frac{w_{j-1}+\rho}{m}\rfloor)(j\nu-m)}&\geq(\displaystyle{\lfloor\frac{w_{m-1}+\rho}{m}\rfloor-\lfloor\frac{w_{1}+\rho}{m}\rfloor)(2\nu-m).}
\end{array}
$  Hence,\newline
$
\begin{array}{lll}
\displaystyle{\sum_{j=1}^{m-1}(\lfloor\frac{w_j+\rho}{m}\rfloor-\lfloor\frac{w_{j-1}+\rho}{m}\rfloor)(j\nu-m)+\rho
}&=\displaystyle{(\lfloor\frac{w_1+\rho}{m}\rfloor-\lfloor\frac{w_0+\rho}{m}\rfloor)(\nu-m)}\\&+\displaystyle{\sum_{j=2}^{m-1}(\lfloor\frac{w_j+\rho}{m}\rfloor-\lfloor\frac{w_{j-1}+\rho}{m}\rfloor)(j\nu-m)+\rho}\\ \\ 
&\geq
\displaystyle{\lfloor\frac{w_1+\rho}{m}\rfloor(\nu-m)}\\&+\displaystyle{(\lfloor\frac{w_{m-1}+\rho}{m}\rfloor-\lfloor\frac{w_{1}+\rho}{m}\rfloor)(2\nu-m)+\rho}\\ \\

&=
\displaystyle{\lfloor\frac{w_1+\rho}{m}\rfloor(\nu-m+(2\nu-m)-(2\nu-m))}\\&+\displaystyle{(\lfloor\frac{w_{m-1}+\rho}{m}\rfloor-\lfloor\frac{w_{1}+\rho}{m}\rfloor)(2\nu-m)+\rho}

\end{array}
$

$
\begin{array}{lll}
\hspace{6.5cm} &=\displaystyle{\lfloor\frac{w_1+\rho}{m}\rfloor(3\nu-2m)}\\ &
+\displaystyle{(\lfloor\frac{w_{m-1}+\rho}{m}\rfloor-\lfloor\frac{w_{1}+\rho}{m}\rfloor-\lfloor\frac{w_{1}+\rho}{m}\rfloor)(2\nu-m)+\rho}\\ \\
&=\displaystyle{\lfloor\frac{w_1+\rho}{m}\rfloor(m-3)}\\ &
+\displaystyle{(\lfloor\frac{w_{m-1}+\rho}{m}\rfloor-\lfloor\frac{w_{1}+\rho}{m}\rfloor-\lfloor\frac{w_{1}+\rho}{m}\rfloor)(m-2)+\rho}.\\ \\
\end{array}
$

\noindent Since $m-\nu=1>0=\frac{1(0)}{2}$, then by Lemma \ref{lem1}, it follows that $w_{m-1}\geq w_1+w_1$. Consequently, by Lemma \ref{x}, we have $
\lfloor\frac{w_{m-1}+\rho}{m}\rfloor\geq \lfloor\frac{w_1+\rho}{m}\rfloor+\lfloor\frac{w_1+\rho}{m}\rfloor-1$.\\ 
\begin{itemize}
\item If $\lfloor\frac{w_{m-1}+\rho}{m}\rfloor- \lfloor\frac{w_1+\rho}{m}\rfloor-\lfloor\frac{w_1+\rho}{m}\rfloor=-1$. Then by Lemma \ref{x}, we have $
\lfloor\frac{w_1+\rho}{m}\rfloor\geq 2$. Since $m\geq 5$, then

$\begin{array}{lll}
\displaystyle{\sum_{j=1}^{m-1}(\lfloor\frac{w_j+\rho}{m}\rfloor-\lfloor\frac{w_{j-1}+\rho}{m}\rfloor)(j\nu-m)+\rho
}
&\geq 2(m-3)-(m-2)+\rho\geq 0.
\end{array}$

\item If $\lfloor\frac{w_{m-1}+\rho}{m}\rfloor- \lfloor\frac{w_1+\rho}{m}\rfloor-\lfloor\frac{w_1+\rho}{m}\rfloor\geq 0$.
\noindent Since $m\geq 5$, then

$\begin{array}{lll}
\displaystyle{\sum_{j=1}^{m-1}(\lfloor\frac{w_j+\rho}{m}\rfloor-\lfloor\frac{w_{j-1}+\rho}{m}\rfloor)(j\nu-m)+\rho
}
&\geq (m-3)+\rho\geq 0.
\end{array}$
\end{itemize}

\noindent Using  Proposition \ref{prop2}, we get that $S$ satisfies Wilf's Conjecture if $m-\nu=1$.  
\item[\textbf{\textit{ii})}] If $m-\nu\in\{2,3\}$. By taking $\alpha=2$ in \eqref{4.1 2}, we get \newline $
\begin{array}{lll}
\displaystyle{\sum_{j={3}}^{m-1}(\lfloor\frac{w_j+\rho}{m}\rfloor-\lfloor\frac{w_{j-1}+\rho}{m}\rfloor)(j\nu-m)}&\geq\displaystyle{(\lfloor\frac{w_{m-1}+\rho}{m}\rfloor-\lfloor\frac{w_{2}+\rho}{m}\rfloor)(3\nu-m).}
\end{array}
$ Hence,
\begin{equation}\label{c2}
\begin{array}{lll}
\displaystyle{\sum_{j=1}^{m-1}(\lfloor\frac{w_j+\rho}{m}\rfloor-\lfloor\frac{w_{j-1}+\rho}{m}\rfloor)(j\nu-m)+\rho
}&=\displaystyle{(\lfloor\frac{w_1+\rho}{m}\rfloor-\lfloor\frac{w_0+\rho}{m}\rfloor)(\nu-m)}\\&\displaystyle{+ (\lfloor\frac{w_2+\rho}{m}\rfloor-\lfloor\frac{w_1+\rho}{m}\rfloor)(2\nu-m)}\\ 
&\displaystyle{+\sum_{j={3}}^{m-1}(\lfloor\frac{w_j+\rho}{m}\rfloor-\lfloor\frac{w_{j-1}+\rho}{m}\rfloor)(j\nu-m)+\rho}\\ \\

&\geq \displaystyle{\lfloor\frac{w_1+\rho}{m}\rfloor(-\nu)+ \lfloor\frac{w_2+\rho}{m}\rfloor(2\nu-m)}\\&+\displaystyle{(\lfloor\frac{w_{m-1}+\rho}{m}\rfloor-\lfloor\frac{w_{2}+\rho}{m}\rfloor)(3\nu-m)+\rho}\\ \\
&= \displaystyle{\lfloor\frac{w_1+\rho}{m}\rfloor(-\nu+(3\nu-m)-(3\nu-m))}\\&\displaystyle{+ \lfloor\frac{w_2+\rho}{m}\rfloor(2\nu-m)}\\&+\displaystyle{(\lfloor\frac{w_{m-1}+\rho}{m}\rfloor-\lfloor\frac{w_{2}+\rho}{m}\rfloor)(3\nu-m)+\rho}\\ \\

&=\displaystyle{\lfloor\frac{w_1+\rho}{m}\rfloor(2\nu-m)}+\displaystyle{\lfloor\frac{w_2+\rho}{m}\rfloor(2\nu-m)}\\ &
+\displaystyle{(\lfloor\frac{w_{m-1}+\rho}{m}\rfloor-\lfloor\frac{w_{1}+\rho}{m}\rfloor-\lfloor\frac{w_{2}+\rho}{m}\rfloor)(3\nu-m)}\\&+\rho.
\end{array}
\end{equation}

\noindent Since $m-\nu\in\{2,3\}>1$, by Lemma \ref{lem1}, we have $w_{m-1}\geq w_1+w_2$. It follows from Lemma \ref{x} that $\lfloor\frac{w_{m-1}+\rho}{m}\rfloor\geq \lfloor\frac{w_1+\rho}{m}\rfloor+\lfloor\frac{w_2+\rho}{m}\rfloor-1$.
\begin{itemize}

\item If $m-\nu=2$. then we may assume that $m=\nu+2\geq 6$. Now \eqref{c2} gives,\newline

$
\begin{array}{lll}
\displaystyle{\sum_{j=1}^{m-1}(\lfloor\frac{w_j+\rho}{m}\rfloor-\lfloor\frac{w_{j-1}+\rho}{m}\rfloor)(j\nu-m)+\rho
}&\geq\displaystyle{\lfloor\frac{w_1+\rho}{m}\rfloor(m-4)}+\displaystyle{\lfloor\frac{w_2+\rho}{m}\rfloor(m-4)}\\&
+\displaystyle{(\lfloor\frac{w_{m-1}+\rho}{m}\rfloor-\lfloor\frac{w_{1}+\rho}{m}\rfloor-\lfloor\frac{w_{2}+\rho}{m}\rfloor)(2m-6)+\rho}.\\ \\
\end{array}
$

\begin{itemize}
\item If $\lfloor\frac{w_{m-1}+\rho}{m}\rfloor- \lfloor\frac{w_1+\rho}{m}\rfloor-\lfloor\frac{w_2+\rho}{m}\rfloor=-1$. Then by Lemma \ref{x} we have, $
\lfloor\frac{w_1+\rho}{m}\rfloor\geq 2\text{  and } \lfloor\frac{w_2+\rho}{m}\rfloor\geq 2$. Since $m\geq 6$, then 

$\begin{array}{lll}
\displaystyle{\sum_{j=1}^{m-1}(\lfloor\frac{w_j+\rho}{m}\rfloor-\lfloor\frac{w_{j-1}+\rho}{m}\rfloor)(j\nu-m)+\rho
}
&\geq 2(m-4)+2(m-4)-(2m-6)+\rho\geq 0.
\end{array}$

\item If $\lfloor\frac{w_{m-1}+\rho}{m}\rfloor- \lfloor\frac{w_1+\rho}{m}\rfloor-\lfloor\frac{w_2+\rho}{m}\rfloor\geq 0$.
\noindent Since $m\geq 6$, then

$\begin{array}{lll}
\displaystyle{\sum_{j=1}^{m-1}(\lfloor\frac{w_j+\rho}{m}\rfloor-\lfloor\frac{w_{j-1}+\rho}{m}\rfloor)(j\nu-m)+\rho
}
&\geq (m-4)+(m-4)+\rho\geq 0.
\end{array}$
\end{itemize}

\noindent Using  Proposition \ref{prop2}, we get that $S$ satisfies Wilf's Conjecture if $m-\nu=2$.
\item If $m-\nu=3$, then we may assume that $m=\nu+3\geq 7$. 
\noindent Now \eqref{c2} gives,\newline
$
\begin{array}{lllll}
\displaystyle{\sum_{j=1}^{m-1}(\lfloor\frac{w_j+\rho}{m}\rfloor-\lfloor\frac{w_{j-1}+\rho}{m}\rfloor)(j\nu-m)+\rho
}&\geq\displaystyle{\lfloor\frac{w_1+\rho}{m}\rfloor(m-6)}+\displaystyle{\lfloor\frac{w_2+\rho}{m}\rfloor(m-6)}
\\ &+\displaystyle{(\lfloor\frac{w_{m-1}+\rho}{m}\rfloor-\lfloor\frac{w_{1}+\rho}{m}\rfloor-\lfloor\frac{w_{2}+\rho}{m}\rfloor)(2m-9)+\rho.}\\ \\
\end{array}
$

\begin{itemize}
\item If $\lfloor\frac{w_{m-1}+\rho}{m}\rfloor- \lfloor\frac{w_1+\rho}{m}\rfloor-\lfloor\frac{w_2+\rho}{m}\rfloor=-1$. Then by Lemma \ref{x} we have, $
\lfloor\frac{w_1+\rho}{m}\rfloor\geq 2, \lfloor\frac{w_2+\rho}{m}\rfloor\geq 2\text{  and } \rho\geq 1$. Since $m\geq 7$, then 

$\begin{array}{lll}
\displaystyle{\sum_{j=1}^{m-1}(\lfloor\frac{w_j+\rho}{m}\rfloor-\lfloor\frac{w_{j-1}+\rho}{m}\rfloor)(j\nu-m)+\rho
}
&\geq 2(m-6)+2(m-6)-(2m-9)+1\geq 0.
\end{array}$

\item If $\lfloor\frac{w_{m-1}+\rho}{m}\rfloor- \lfloor\frac{w_1+\rho}{m}\rfloor-\lfloor\frac{w_2+\rho}{m}\rfloor\geq 0$.
\noindent Since $m\geq 7$, then

$\begin{array}{lll}
\displaystyle{\sum_{j=1}^{m-1}(\lfloor\frac{w_j+\rho}{m}\rfloor-\lfloor\frac{w_{j-1}+\rho}{m}\rfloor)(j\nu-m)+\rho
}
&\geq (m-6)+(m-6)+\rho\geq 0.
\end{array}$
\end{itemize}

\noindent Using Proposition \ref{prop2}, we get that $S$ satisfies Wilf's Conjecture if $m-\nu=3$.  
  \end{itemize}
\end{enumerate}
\noindent Thus Wilf's Conjecture holds if $m-\nu\leq 3$. $ \quad \square $
\end{demostracion}

\noindent The next Corollary covers the result of Sammartano for numerical semigroups with $2\nu\geq m$ (\cite{Sammartano1}) using Corollary \ref{cor1} and Lemma \ref{thm1}.
\begin{corolario}\label{cor2} {\rm Let $S$ be a numerical semigroup with multiplicity $m$ and embedding dimension $\nu$.
If $2\nu\geq m$, then $S$ satisfies Wilf's Conjecture.
}
\end{corolario}
\begin{demostracion}{.} If $m-\nu>3$ and $2\nu\geq m$, then by Corollary \ref{cor1} Wilf's Conjecture holds. If $m-\nu\leq 3$, by Lemma \ref{thm1}, $S$ satisfies Wilf's Conjecture. $\quad \square$
\end{demostracion}
\medskip
\noindent In the following Corollary we will deduce Wilf's Conjecture for numerical semigroups with $m-\nu=4$. This will enable us later to prove the conjecture for those with $(2+\frac{1}{q})\nu\geq m$. 
\begin{corolario}\label{corx} {\rm Let $S$ be a numerical semigroup with multiplicity $m$ and embedding dimension $\nu$. If $m-\nu=4$, then $S$ satisfies Wilf's Conjecture.
}
\end{corolario}
\begin{demostracion}{.} Since Wilf's conjecture holds for $\nu\leq 3$ (\cite{Dobbs1}), then we may assume that $\nu\geq 4$. Hence, $\nu\geq m-\nu$. Consequently, $2\nu\geq m$. Hence, $S$ satisfies Wilf's Conjecture. $\quad \square$
\end{demostracion}
\medskip
\noindent The following technical Lemma will be used through the paper.

\begin{lema}\label{lem2} {\rm Let $S$ be a numerical semigroup with multiplicity $m$ and embedding dimension $\nu$. Let $w_0=0<w_1<\ldots<w_{m-1}$ be the elements of Ap$(S,m)$. If $m-\nu\geq (^\alpha_2)-1=\frac{\alpha(\alpha-1)}{2}-1$ for some $3\leq\alpha\leq m-2$, then $w_{m-1}\geq w_{1}+w_{\alpha}$ or $w_{m-1}\geq w_{\alpha-2}+w_{\alpha-1}$. }
\end{lema}

\medskip

\begin{demostracion}{.} Suppose by the way of contradiction that $w_{m-1}< w_{1}+w_{\alpha}$ and $w_{m-1}< w_{\alpha-2}+w_{\alpha-1}$. Let\newline $w\in$Ap$(S,m)\setminus$min(Ap$(S,m))$, then $w\leq w_{m-1}$ and $w=w_i+w_j$ for some $w_i,w_j\in$Ap$(S,m)\setminus\{0\}$. In this case, the only possible values of $w$ are $\{w_i+w_j;1\leq i\leq j\leq\alpha-1\}\setminus\{w_{\alpha-2}+w_{\alpha-1},w_{\alpha-1}+w_{\alpha-1}\}.$ Consequently, $m-\nu=\vert$Ap$(S,m)\setminus$min(Ap$(S,m))\vert\leq\frac{\alpha(\alpha-1)}{2}-2$. But $\frac{\alpha(\alpha-1)}{2}-2<\frac{\alpha(\alpha-1)}{2}-1$, which contradicts the hypothesis. Hence, $w_{m-1}\geq w_{1}+w_{\alpha}$ or $w_{m-1}\geq w_{\alpha-2}+w_{\alpha-1}$. $\quad \square$

\end{demostracion}

\medskip
\noindent In the next theorem, we will show that Wilf's Conjecture holds for numerical semigroups with $m-\nu=5$. 

\begin{teorema}\label{thm7} {\rm Let $S$ be a numerical semigroup with multiplicity $m$ and embedding dimension $\nu$. If $m-\nu=5$, then $S$ satisfies Wilf's Conjecture.}
\end{teorema}

\medskip

\begin{demostracion}{.} Let $m-\nu=5$. Since Wilf's Conjecture holds for $2\nu\geq m$, then we may assume that $2\nu<m$. This implies that $\nu<5$. Since the case $\nu\leq 3$ is known (\cite{Dobbs1}), then we shall assume that $\nu=4$. This also implies that $m=9$.

\noindent Since $m-\nu=5=\frac{4(3)}{2}-1$, by Lemma \ref{lem2}, it follows that $w_{8}\geq w_2+w_3$ or $w_{8}\geq w_1+w_4$.
 
\begin{enumerate}
 
\item[\textbf{\textit{i})}] If $w_{8}\geq w_2+w_3$. By taking $\alpha=3$ in \eqref{4.1 2} ($m=9$, $\nu=4$), we get

$
\begin{array}{lll}
\displaystyle{\sum_{j={4}}^{8}(\lfloor\frac{w_j+\rho}{9}\rfloor-\lfloor\frac{w_{j-1}+\rho}{9}\rfloor)(4j-9)}&\geq(\displaystyle{\lfloor\frac{w_{8}+\rho}{9}\rfloor-\lfloor\frac{w_{3}+\rho}{9}\rfloor)(16-9)}&=(\displaystyle{\lfloor\frac{w_{8}+\rho}{9}\rfloor-\lfloor\frac{w_{3}+\rho}{9}\rfloor)(7)}.
\end{array}
$

\noindent Hence,
\begin{equation}\label{a1}
\begin{array}{lllll}
\displaystyle{\sum_{j=1}^{8}(\lfloor\frac{w_j+\rho}{9}\rfloor-\lfloor\frac{w_{j-1}+\rho}{9}\rfloor)(4j-9)+\rho}
&=\displaystyle{(\lfloor\frac{w_1+\rho}{9}\rfloor-\lfloor\frac{w_0+\rho}{9}\rfloor)(-5) }\\&\displaystyle{+ (\lfloor\frac{w_2+\rho}{9}\rfloor-\lfloor\frac{w_1+\rho}{9}\rfloor)(-1)}\\&\displaystyle{+(\lfloor\frac{w_3+\rho}{9}\rfloor-\lfloor\frac{w_2+\rho}{9}\rfloor)(3)}\\&\displaystyle{+\sum_{j={4}}^{8}(\lfloor\frac{w_j+\rho}{9}\rfloor-\lfloor\frac{w_{j-1}+\rho}{9}\rfloor)(4j-9)}+\rho
\\ \\
&\geq\displaystyle{\lfloor\frac{w_1+\rho}{9}\rfloor(-4) + \lfloor\frac{w_2+\rho}{9}\rfloor(-4)}\displaystyle{+\lfloor\frac{w_3+\rho}{9}\rfloor(3)}\\&\displaystyle{+(\lfloor\frac{w_{8}+\rho}{9}\rfloor-\lfloor\frac{w_{3}+\rho}{9}\rfloor)(7)+\rho}
\\ \\
&\geq\displaystyle{\bigg(\lfloor\frac{w_2+\rho}{m}\rfloor((\frac{-3}{4})4)+\lfloor\frac{w_3+\rho}{9}\rfloor((\frac{-1}{4})4)\bigg) }\\&\displaystyle{+ \lfloor\frac{w_2+\rho}{9}\rfloor(-4)}\displaystyle{+\lfloor\frac{w_3+\rho}{9}\rfloor(3)}\\&\displaystyle{+(\lfloor\frac{w_{8}+\rho}{9}\rfloor-\lfloor\frac{w_{3}+\rho}{9}\rfloor)(7)+\rho}

\end{array}
\end{equation}

\begin{equation*}
\begin{array}{lllll} 

\hspace{7.5cm}&= \displaystyle{\lfloor\frac{w_2+\rho}{9}\rfloor(-7)}+\displaystyle{\lfloor\frac{w_3+\rho}{9}\rfloor(2)}\\&+
\displaystyle{(\lfloor\frac{w_{8}+\rho}{9}\rfloor-\lfloor\frac{w_{3}+\rho}{9}\rfloor)(7)+\rho}\\ \\
&=\displaystyle{\lfloor\frac{w_3+\rho}{9}\rfloor(2)}\\&+\displaystyle{(\lfloor\frac{w_{8}+\rho}{9}\rfloor-\lfloor\frac{w_{2}+\rho}{9}\rfloor-\lfloor\frac{w_{3}+\rho}{9}\rfloor)(7)+\rho.}\\ \\
\end{array}
\end{equation*}
\noindent Since $w_{8}\geq w_2+w_3$, by Lemma \ref{x}, it follows that $\lfloor\frac{w_{8}+\rho}{9}\rfloor\geq \lfloor\frac{w_2+\rho}{9}\rfloor+\lfloor\frac{w_3+\rho}{9}\rfloor-1$.

\begin{itemize}
\item If $\lfloor\frac{w_{8}+\rho}{9}\rfloor- \lfloor\frac{w_2+\rho}{9}\rfloor-\lfloor\frac{w_3+\rho}{9}\rfloor\geq 0$,  then \eqref{a1} gives $\displaystyle{\sum_{j=1}^{8}(\lfloor\frac{w_j+\rho}{9}\rfloor-\lfloor\frac{w_{j-1}+\rho}{9}\rfloor)(4j-9)+\rho}\geq 0.$

\item If $\lfloor\frac{w_{8}+\rho}{9}\rfloor- \lfloor\frac{w_2+\rho}{9}\rfloor-\lfloor\frac{w_3+\rho}{9}\rfloor=-1$. By Lemma \ref{x}, we have $\rho\geq 1$. 

\noindent Since for $q\leq 3$ Wilf's Conjecture is solved (\cite{Eliahou}, \cite{Kaplan1}), then may assume that $q\geq4$. Since $\lfloor\frac{w_2+\rho}{9}\rfloor\leq\lfloor\frac{w_3+\rho}{9}\rfloor$ and $\lfloor\frac{w_2+\rho}{9}\rfloor+\lfloor\frac{w_3+\rho}{9}\rfloor=\lfloor\frac{w_{8}+\rho}{9}\rfloor+1=q+1$, in this case it follows that $\lfloor\frac{w_3+\rho}{9}\rfloor\geq 3.$

\noindent Now \eqref{a1} gives, $
\displaystyle{\sum_{j=1}^{8}(\lfloor\frac{w_j+\rho}{9}\rfloor-\lfloor\frac{w_{j-1}+\rho}{9}\rfloor)(4j-m)+\rho}
\geq 3(2)-7+1\geq 0.$

\end{itemize}

\noindent Using  Proposition \ref{prop2}, we get that $S$ satisfies Wilf's Conjecture in this case.

\item[\textbf{\textit{ii})}] If $w_{8}\geq w_1+w_4$. We may assume that $w_{8}< w_2+w_3$, since otherwise we are back to case 1. Hence, the possible values of $w\in$ Ap$(S,9)\setminus$min(Ap($S,9))$ are $\{w_1+w_j; 1\leq j\leq 7\}\cup\{w_2+w_2\}$. 

\begin{itemize}

\item Recall that an element $x$ of the Ap\'ery set of $S$ belongs to max(Ap($S,m))$ if and only if $w_i\neq x+w_j$ for all $w_i,w_j\in$Ap($S,m)\setminus\{0\}$. If Ap$(S,9)\setminus$min(Ap($S,9))\subseteq \{w_1+w_j; 1\leq j\leq 7\}$, then there exists at least five elements in Ap$(S,9)$ that are not maximal, hence $t(S)=\vert\{$max(Ap($S,9))-9\}\vert\leq3=\nu-1$. Consequently, $S$ satisfies Wilf's Conjecture (\cite{Dobbs1} Proposition 2.3).

\item If $w_2+w_2\in$ Ap$(S,9)\setminus$min(Ap($S,9))$, then $w_{8}\geq w_2+w_2$. By Lemma \ref{x} we have $\lfloor\frac{w_{8}+\rho}{9}\rfloor\geq 2\lfloor\frac{w_2+\rho}{9}\rfloor-1$. In particular, \begin{equation}\label{k}
\lfloor\frac{w_2+\rho}{9}\rfloor\leq\frac{q+1}{2}.\end{equation}\newline
\noindent By taking $\alpha=4$ in \eqref{4.1 2} ($m=9,$ $\nu=4$), we get \newline
$
\begin{array}{llll}
\displaystyle{\sum_{j={5}}^{8}(\lfloor\frac{w_j+\rho}{9}\rfloor-\lfloor\frac{w_{j-1}+\rho}{9}\rfloor)(4j-9)}&\geq(\displaystyle{\lfloor\frac{w_{8}+\rho}{9}\rfloor-\lfloor\frac{w_{4}+\rho}{9}\rfloor)(11).}
\end{array}
$ Now using  \eqref{k}, we get

\begin{equation}\label{a2}
\begin{array}{lllll}
\displaystyle{\sum_{j=1}^{8}(\lfloor\frac{w_j+\rho}{9}\rfloor-\lfloor\frac{w_{j-1}+\rho}{9}\rfloor)(4j-9)+\rho}
&=\displaystyle{(\lfloor\frac{w_1+\rho}{9}\rfloor-\lfloor\frac{w_0+\rho}{9}\rfloor)(-5)}\\&\displaystyle{ + (\lfloor\frac{w_2+\rho}{9}\rfloor-\lfloor\frac{w_1+\rho}{9}\rfloor)(-1)}\\&\displaystyle{+(\lfloor\frac{w_3+\rho}{9}\rfloor-\lfloor\frac{w_2+\rho}{m}\rfloor)(3)}\\&\displaystyle{ + (\lfloor\frac{w_4+\rho}{9}\rfloor-\lfloor\frac{w_3+\rho}{9}\rfloor)(7)}\\&\displaystyle{+\sum_{j={5}}^{8}(\lfloor\frac{w_j+\rho}{9}\rfloor-\lfloor\frac{w_{j-1}+\rho}{9}\rfloor)(4j-9)+\rho}
\end{array}
\end{equation}

\begin{equation*}
\begin{array}{lllll}
\hspace{7.5cm}&\geq\displaystyle{\lfloor\frac{w_1+\rho}{9}\rfloor(-4)+\lfloor\frac{w_2+\rho}{9}\rfloor(-4)+\lfloor\frac{w_3+\rho}{9}\rfloor(-4)}\\ &\displaystyle{ + \lfloor\frac{w_4+\rho}{9}\rfloor(7)}\displaystyle{+(\lfloor\frac{w_{8}+\rho}{9}\rfloor-\lfloor\frac{w_{4}+\rho}{9}\rfloor)(11)}+\rho
\\ \\
&\geq\displaystyle{\lfloor\frac{w_1+\rho}{9}\rfloor(-4) + (\frac{q+1}{2})(-4)+\lfloor\frac{w_4+\rho}{9}\rfloor(-4)}\\ &\displaystyle{ + \lfloor\frac{w_4+\rho}{9}\rfloor(7)}\displaystyle{+(\lfloor\frac{w_{8}+\rho}{9}\rfloor-\lfloor\frac{w_{4}+\rho}{9}\rfloor)(11)}+\rho
\\ \\
&=\displaystyle{\lfloor\frac{w_1+\rho}{9}\rfloor(-4) -2 (q+1)+\lfloor\frac{w_4+\rho}{9}\rfloor(3)}\\ &\displaystyle{+(\lfloor\frac{w_{8}+\rho}{9}\rfloor-\lfloor\frac{w_{4}+\rho}{9}\rfloor)(11)+\rho}
\\ \\
&=\displaystyle{\lfloor\frac{w_1+\rho}{9}\rfloor(-4+11-11) -2(q+1)}\\ &\displaystyle{+\lfloor\frac{w_4+\rho}{9}\rfloor(3)+(\lfloor\frac{w_{8}+\rho}{9}\rfloor-\lfloor\frac{w_{4}+\rho}{9}\rfloor)(11)+\rho}
\\ \\
&=\displaystyle{\lfloor\frac{w_1+\rho}{9}\rfloor(7) -2(q+1)+\lfloor\frac{w_4+\rho}{9}\rfloor(3)}\\ &\displaystyle{+(\lfloor\frac{w_{8}+\rho}{9}\rfloor-\lfloor\frac{w_{4}+\rho}{9}\rfloor-\lfloor\frac{w_1+\rho}{9}\rfloor)(11)+\rho}
\\ \\

&=\displaystyle{(\lfloor\frac{w_1+\rho}{9}\rfloor+\lfloor\frac{w_4+\rho}{9}\rfloor)(3)}+\displaystyle{\lfloor\frac{w_1+\rho}{9}\rfloor(4)}-\displaystyle{2(q+1)
}\\&+\displaystyle{(\lfloor\frac{w_{8}+\rho}{9}\rfloor-\lfloor\frac{w_{1}+\rho}{9}\rfloor-\lfloor\frac{w_{4}+\rho}{9}\rfloor)(11)}+\rho.\\ \\
\end{array}
\end{equation*}

\noindent We have $w_{8}\geq w_1+w_4$, then by Lemma \ref{x} $\lfloor\frac{w_{8}+\rho}{9}\rfloor\geq \lfloor\frac{w_1+\rho}{9}\rfloor+\lfloor\frac{w_4+\rho}{9}\rfloor-1$.
\begin{itemize}

\item If $\lfloor\frac{w_{8}+\rho}{9}\rfloor- \lfloor\frac{w_1+\rho}{9}\rfloor-\lfloor\frac{w_4+\rho}{9}\rfloor\geq 0$. Let $x=\lfloor\frac{w_{8}+\rho}{9}\rfloor- \lfloor\frac{w_1+\rho}{9}\rfloor-\lfloor\frac{w_4+\rho}{9}\rfloor$. Hence, $x\geq0$ and $\lfloor\frac{w_1+\rho}{9}\rfloor+\lfloor\frac{w_4+\rho}{9}\rfloor=q-x$. Then \eqref{a2} gives,

$
\begin{array}{lll}
\displaystyle{\sum_{j=1}^{8}(\lfloor\frac{w_j+\rho}{9}\rfloor-\lfloor\frac{w_{j-1}+\rho}{9}\rfloor)(4j-9)+\rho}
&\geq \displaystyle{(q-x)(3)}+\displaystyle{4}-\displaystyle{2(q+1)
+11x}+\rho=q+8x+2+\rho
\geq 0.
\end{array}
$

\item If $\lfloor\frac{w_{8}+\rho}{9}\rfloor- \lfloor\frac{w_1+\rho}{9}\rfloor-\lfloor\frac{w_4+\rho}{9}\rfloor=-1$. Then $\lfloor\frac{w_1+\rho}{m}\rfloor+\lfloor\frac{w_4+\rho}{9}\rfloor=q+1$. By Lemma \ref{x}, we have $
\lfloor\frac{w_1+\rho}{9}\rfloor\geq 2$  and $\rho\geq 1.$ Since $q\geq 1$, then \eqref{a2} gives,\newline
$
\begin{array}{lll}
\displaystyle{\sum_{j=1}^{8}(\lfloor\frac{w_j+\rho}{9}\rfloor-\lfloor\frac{w_{j-1}+\rho}{9}\rfloor)(4j-9)+\rho}
&\geq \displaystyle{(q+1)(3)}+\displaystyle{8}-\displaystyle{2(q+1)
-11}+1= q-1\geq 0.
\end{array}
$

\end{itemize}

\noindent Using  Proposition \ref{prop2}, we get that $S$ satisfies Wilf's Conjecture in this case.  

\end{itemize}

\noindent Thus, Wilf's Conjecture holds if $m-\nu=5$. $ \quad \square $

\end{enumerate}

\end{demostracion}

\medskip
\noindent In the next corollary, we will deduce the conjecture for $m=9$.
\begin{corolario} {\rm If $S$ is a numerical Semigroup with multiplicity $m=9$, then $S$ satisfies Wilf's Conjecture.}
\end{corolario}
\medskip
\begin{demostracion}{.} By Lemma \ref{thm1}, Corollary \ref{corx} and Theorem \ref{thm7}, we may assume that $m-\nu>5$, hence $\nu<m-5=4$. By (\cite{Dobbs1}) $S$ satisfies Wilf's Conjecture. $\quad \square$

\end{demostracion}

\medskip
\noindent The following Lemma will enable us later to show that Wilf's Conjecture holds for numerical semigroups with $(2+\frac{1}{q})\nu\geq m$. 
\begin{lema}\label{thm8}{\rm Let $S$ be a numerical Semigroup with multiplicity $m$, embedding dimension $\nu$ and conductor $f+1=qm-\rho$ for some $q\in \N,$ $0\leq \rho\leq m-1$. If $m-\nu=6$ and $(2+\frac{1}{q})\nu\geq m$, then $S$ satisfies Wilf's Conjecture.
}
\end{lema}

\medskip

\begin{demostracion}{.}  Since $m-\nu=6\geq\frac{4(3)}{2}-1$, by Lemma \ref{lem2}, it follows that $w_{m-1}\geq w_1+w_4$ or $w_{m-1}\geq w_2+w_3$.
 
\begin{enumerate}
 
\item[\textbf{\textit{i})}] If $w_{m-1}\geq w_1+w_4$. By hypothesis $(2+\frac{1}{q})\nu\geq m$ and Theorem \ref{thm2} Wilf's Conjecture holds in this case.

\item[\textbf{\textit{ii})}] If $w_{m-1}\geq w_2+w_3$. We may assume that $w_{m-1}< w_1+w_4$, since otherwise we are back to case 1. Hence, Ap$(S,m)\setminus$min(Ap($S,m))=\{w_1+w_1,w_1+w_2,w_1+w_3,w_2+w_2,w_2+w_3,w_3+w_3\}$. 

\noindent By taking $\alpha=3$ in \eqref{4.1 2}, we get
$
\begin{array}{llll}
\displaystyle{\sum_{j={4}}^{m-1}(\lfloor\frac{w_j+\rho}{m}\rfloor-\lfloor\frac{w_{j-1}+\rho}{m}\rfloor)(j\nu-m)}&\geq(\displaystyle{\lfloor\frac{w_{m-1}+\rho}{m}\rfloor-\lfloor\frac{w_{3}+\rho}{m}\rfloor)(4\nu-m).}
\end{array}
$
\noindent Hence,

\begin{equation}\label{b1}
\begin{array}{lllll}
\displaystyle{\sum_{j=1}^{m-1}(\lfloor\frac{w_j+\rho}{m}\rfloor-\lfloor\frac{w_{j-1}+\rho}{m}\rfloor)(j\nu-m)+\rho}
&=\displaystyle{(\lfloor\frac{w_1+\rho}{m}\rfloor-\lfloor\frac{w_0+\rho}{m}\rfloor)(\nu-m)}\\ &\displaystyle{ + (\lfloor\frac{w_2+\rho}{m}\rfloor-\lfloor\frac{w_1+\rho}{m}\rfloor)(2\nu-m)}\\&\displaystyle{+(\lfloor\frac{w_3+\rho}{m}\rfloor-\lfloor\frac{w_2+\rho}{m}\rfloor)(3\nu-m)}\\ &\displaystyle{+\sum_{j={4}}^{m-1}(\lfloor\frac{w_j+\rho}{m}\rfloor-\lfloor\frac{w_{j-1}+\rho}{m}\rfloor)(j\nu-m)+\rho}
\\ \\
&\geq\displaystyle{\lfloor\frac{w_1+\rho}{m}\rfloor(-\nu)}\displaystyle{ + \lfloor\frac{w_2+\rho}{m}\rfloor(-\nu)}\\&\displaystyle{+\lfloor\frac{w_3+\rho}{m}\rfloor(3\nu-m)}\\ &\displaystyle{+(\lfloor\frac{w_{m-1}+\rho}{m}\rfloor-\lfloor\frac{w_{3}+\rho}{m}\rfloor)(4\nu-m)+\rho}
\\ \\
&\geq\displaystyle{\big(\lfloor\frac{w_2+\rho}{m}\rfloor(\frac{-\nu}{2})+\lfloor\frac{w_3+\rho}{m}\rfloor(\frac{-\nu}{2})\big)}\\&\displaystyle{ + \lfloor\frac{w_2+\rho}{m}\rfloor(-\nu)}\displaystyle{+\lfloor\frac{w_3+\rho}{m}\rfloor(3\nu-m)}\\& \displaystyle{+(\lfloor\frac{w_{m-1}+\rho}{m}\rfloor-\lfloor\frac{w_{3}+\rho}{m}\rfloor)(4\nu-m)}+\rho
\\ \\
&=\displaystyle{\lfloor\frac{w_2+\rho}{m}\rfloor(\frac{-3\nu}{2})+\lfloor\frac{w_3+\rho}{m}\rfloor(\frac{5\nu}{2}-m)}\\&\displaystyle{+(\lfloor\frac{w_{m-1}+\rho}{m}\rfloor-\lfloor\frac{w_{3}+\rho}{m}\rfloor)(4\nu-m)+\rho}
\\ \\
&=\displaystyle{\lfloor\frac{w_2+\rho}{m}\rfloor(\frac{-3\nu}{2}+(4\nu-m)-(4\nu-m))}\\&\displaystyle{+\lfloor\frac{w_3+\rho}{m}\rfloor(\frac{5\nu}{2}-m)}\\&\displaystyle{+(\lfloor\frac{w_{m-1}+\rho}{m}\rfloor-\lfloor\frac{w_{3}+\rho}{m}\rfloor)(4\nu-m)+\rho}
\\ \\

&=\displaystyle{\lfloor\frac{w_2+\rho}{m}\rfloor(\frac{5\nu}{2}-m)}+\displaystyle{\lfloor\frac{w_3+\rho}{m}\rfloor(\frac{5\nu}{2}-m)}\\
&+\displaystyle{(\lfloor\frac{w_{m-1}+\rho}{m}\rfloor-\lfloor\frac{w_{2}+\rho}{m}\rfloor-\lfloor\frac{w_{3}+\rho}{m}\rfloor)(4\nu-m)}\\&+\rho\\ \\

\end{array}
\end{equation}

$
\begin{array}{lllll}

\hspace{7.5cm}&=\displaystyle{\lfloor\frac{w_2+\rho}{m}\rfloor(\frac{3\nu}{2}-6)}+\displaystyle{\lfloor\frac{w_3+\rho}{m}\rfloor(\frac{3\nu}{2}-6)}\\
&+\displaystyle{(\lfloor\frac{w_{m-1}+\rho}{m}\rfloor-\lfloor\frac{w_{2}+\rho}{m}\rfloor-\lfloor\frac{w_{3}+\rho}{m}\rfloor)(3\nu-6)}\\&+\rho.\\ \\
\end{array}
$

\noindent We have $w_{m-1}\geq w_2+w_3$, by Lemma \ref{x}, it follows that $\lfloor\frac{w_{m-1}+\rho}{m}\rfloor\geq \lfloor\frac{w_2+\rho}{m}\rfloor+\lfloor\frac{w_3+\rho}{m}\rfloor-1$.
\begin{itemize}

\item If $\lfloor\frac{w_{m-1}+\rho}{m}\rfloor- \lfloor\frac{w_2+\rho}{m}\rfloor-\lfloor\frac{w_3+\rho}{m}\rfloor\geq 0$,  using $\nu\geq 4$ in \eqref{b1}
($\nu\leq 3$ is solved \cite{Dobbs1}), we get 

$
\displaystyle{\sum_{j=1}^{m-1}(\lfloor\frac{w_j+\rho}{m}\rfloor-\lfloor\frac{w_{j-1}+\rho}{m}\rfloor)(j\nu-m)+\rho}
\geq 0.
$
\item If $\lfloor\frac{w_{m-1}+\rho}{m}\rfloor- \lfloor\frac{w_2+\rho}{m}\rfloor-\lfloor\frac{w_3+\rho}{m}\rfloor=-1$. Then,\begin{equation}\label{aa3}
\lfloor\frac{w_2+\rho}{m}\rfloor+\lfloor\frac{w_3+\rho}{m}\rfloor=q+1.
\end{equation}
\noindent We have $w_3+w_3\in$Ap$(S,m)\setminus$min(Ap($S,m))$, then $w_{m-1}\geq w_3+w_3$. By Lemma \ref{x}, we have $\lfloor\frac{w_{m-1}+\rho}{m}\rfloor\geq 2\lfloor\frac{w_3+\rho}{m}\rfloor-1$. In particular,  \begin{equation}\label{o}
\lfloor\frac{w_3+\rho}{m}\rfloor\leq\frac{q+1}{2}.\end{equation}

\noindent Since Wilf's Conjecture holds for $q\leq 3$ (\cite{Eliahou}, \cite{Kaplan1}), so we may assume that $q\geq 4$. Since $\lfloor\frac{w_2+\rho}{m}\rfloor\leq\lfloor\frac{w_3+\rho}{m}\rfloor$, by \eqref{aa3} and \eqref{o}, it follows that $\lfloor\frac{w_2+\rho}{m}\rfloor=\lfloor\frac{w_3+\rho}{m}\rfloor=\frac{q+1}{2}$, in particular $q$ is odd, so we have to assume that $ q\geq 5$. Now using Now using \eqref{aa3}, $q\geq 5$ and the hypothesis $(2+\frac{1}{q})\nu\geq m=\nu+6$ in \eqref{b1}, we get

\begin{equation*}
\begin{array}{lll}
\displaystyle{\sum_{j=1}^{m-1}(\lfloor\frac{w_j+\rho}{m}\rfloor-\lfloor\frac{w_{j-1}+\rho}{m}\rfloor)(j\nu-m)+\rho}&\geq\displaystyle{(\lfloor\frac{w_2+\rho}{m}\rfloor+\lfloor\frac{w_3+\rho}{m}\rfloor)(\frac{3\nu}{2}-6)}\\
&+\displaystyle{(\lfloor\frac{w_{m-1}+\rho}{m}\rfloor-\lfloor\frac{w_{2}+\rho}{m}\rfloor-\lfloor\frac{w_{3}+\rho}{m}\rfloor)(3\nu-6)+\rho}\\ \\
&=\displaystyle{(q+1)(\frac{3\nu}{2}-6)-(3\nu-6)+\rho}\\ \\
&=\displaystyle{\nu(\frac{3q}{2}+\frac{3}{2}-3)-6q+\rho}\\ \\
&\geq\displaystyle{\nu(\frac{3q}{2}-\frac{3}{2})-q\nu-\nu+\rho}\hspace{2cm} (6q\leq q\nu+\nu)\\ \\
&=\displaystyle{\nu(\frac{q}{2}-\frac{5}{2})+\rho}
\geq 0.
\end{array}
\end{equation*}
\noindent Using  Proposition \ref{prop2}, we get that $S$ satisfies Wilf's Conjecture in this case. 
\end{itemize}

\noindent Thus, Wilf's Conjecture holds if $m-\nu=6$ and $(2+\frac{1}{q})\nu\geq m$. $ \quad \square $

\end{enumerate}
\end{demostracion}

\medskip
\noindent Next we will generalize a result for Sammartano (\cite{Sammartano1}) and show that Wilf's Conjecture holds for numerical semigroups satisfying $(2+\frac{1}{q})\nu\geq m$, using Lemma \ref{thm1},  Corollary \ref{corx},  Theorem \ref{thm7}, Lemma \ref{thm8} and Corollary \ref{cor1}.
\begin{teorema}\label{cor8} {\rm Let $S$ be a numerical semigroup with multiplicity $m$, embedding dimension $\nu$ and conductor $f+1=qm-\rho$ for some $q\in\N,$ $0\leq \rho\leq m-1$. If $(2+\frac{1}{q})\nu\geq m$, then $S$ satisfies Wilf's Conjecture.
}
\end{teorema}

\begin{demostracion}{.} 
\begin{itemize}
\item If $m-\nu\leq 3$, then by Lemma \ref{thm1} Wilf's Conjecture holds.
\item If $m-\nu=4$, then by Corollary \ref{corx} Wilf's Conjecture holds.
\item If $m-\nu=5$, then by Theorem \ref{thm7} Wilf's Conjecture holds.
\item If $m-\nu=6$ and $(2+\frac{1}{q})\nu\geq m$, then by Lemma \ref{thm8} Wilf's Conjecture holds.
\item If $m-\nu>6$ and $(2+\frac{1}{q})\nu\geq m$, then by Corollary \ref{cor1} Wilf's Conjecture holds. $\quad \square$

\end{itemize}

\end{demostracion}

\begin{exemple} {\rm Consider the following numerical semigroup $S=<13,15,17,19,21,27>$. Note that $2\nu<m$. We have $(2+\frac{1}{q})\nu=(2+\frac{1}{4})6\geq 13=m$. Thus the conditions of Theorem \ref{cor8} are valid, so $S$ satisfies Wilf's Conjecture.}

\end{exemple}

\begin{corolario}
{\rm Let $S$ be a numerical semigroup with multiplicity $m$ and conductor $f+1=qm-\rho$ for some $q\in\N$, $0\leq \rho\leq m-1$. If $m\leq 8+\frac{4}{q}$, then $S$ satisfies Wilf's Conjecture.}
\end{corolario}
\begin{demostracion}{.} If $\nu<4$, then $S$ satisfies Wilf's Conjecture (\cite{Dobbs1}). Hence, we can suppose that $\nu\geq 4$. Thus, $(2+\frac{1}{q})\nu\geq(2+\frac{1}{q})4\geq m$.
 By using Theorem \ref{cor8} $S$ satisfies Wilf's Conjecture . $\quad \square$
\end{demostracion}

\section{Numerical semigroups with $w_{m-1}\geq w_{\alpha-1}+w_{\alpha}$ and $(\frac{\alpha+3}{3})\nu\geq m$} \label{sec5 N.S}
In this section, we will show that if $S$ is a numerical Semigroup such that $w_{m-1}\geq w_{\alpha-1}+w_{\alpha}$ and $(\frac{\alpha+3}{3})\nu\geq m$, then $S$ satisfies Wilf's Conjecture.
\begin{theorema}\label{thm3} {\rm Let $S$ be a numerical semigroup with multiplicity $m$ and embedding dimension $\nu$. Let $w_0=0<w_1<w_2<\ldots<w_{m-1}$ be the elements of Ap($S,m)$. Suppose that $w_{m-1}\geq w_{\alpha-1}+w_{\alpha}$ for some $1<\alpha<m-1$.
If $(\frac{\alpha+3}{3})\nu\geq m$, then $S$ satisfies Wilf's Conjecture.
}
\end{theorema}
\begin{demostracion}{.} We may assume that $\rho\geq \frac{(3-q)\alpha m}{2\alpha+6}$. Indeed, if $0\leq\rho< \frac{(3-q)\alpha m}{2\alpha+6}$, then $q<3$ and Wilf's conjecture holds for this case (\cite{Kaplan1}). We are going to show that $S$ satisfies Wilf's Conjecture by means of Proposition \ref{prop2}. We have,\newline

$
\begin{array}{lll}
\displaystyle{\sum_{j={1}}^{\alpha}(\lfloor\frac{w_j+\rho}{m}\rfloor-\lfloor\frac{w_{j-1}+\rho}{m}\rfloor)}(j\nu-m)&=\displaystyle{\sum_{j={1}}^{\alpha}\lfloor\frac{w_j+\rho}{m}\rfloor(j\nu-m)-\sum_{j=1}^{\alpha}\lfloor\frac{w_{j-1}+\rho}{m}\rfloor}(j\nu-m)\\&=\displaystyle{\sum_{j={1}}^{\alpha}\lfloor\frac{w_j+\rho}{m}\rfloor(j\nu-m)-\sum_{j=0}^{\alpha-1}\lfloor\frac{w_j+\rho}{m}\rfloor}((j+1)\nu-m)\\&= \displaystyle{\lfloor\frac{w_{\alpha}+\rho}{m}\rfloor(\alpha\nu-m)-\lfloor\frac{w_{0}+\rho}{m}\rfloor(\nu-m)-\sum_{j=1}^{\alpha-1}\lfloor\frac{w_{j}+\rho}{m}\rfloor\nu}\\ &= \displaystyle{\lfloor\frac{w_{\alpha}+\rho}{m}\rfloor(\alpha\nu-m)-\lfloor\frac{w_{\alpha-1}+\rho}{m}\rfloor\nu-\sum_{j=1}^{\alpha-2}\lfloor\frac{w_{j}+\rho}{m}\rfloor\nu}
\\& \geq \displaystyle{\lfloor\frac{w_{\alpha}+\rho}{m}\rfloor(\alpha\nu-m)-\lfloor\frac{w_{\alpha-1}+\rho}{m}\rfloor\nu-\sum_{j=1}^{\alpha-2}\bigg(\frac{\lfloor\frac{w_{\alpha}+\rho}{m}\rfloor+\lfloor\frac{w_{\alpha-1}+\rho}{m}\rfloor}{2}\bigg)\nu}
\\& = \displaystyle{\lfloor\frac{w_{\alpha}+\rho}{m}\rfloor(\alpha\nu-m)-\lfloor\frac{w_{\alpha-1}+\rho}{m}\rfloor\nu-(\lfloor\frac{w_{\alpha}+\rho}{m}\rfloor+\lfloor\frac{w_{\alpha-1}+\rho}{m}\rfloor)\frac{(\alpha-2)\nu}{2}}\\ &= \displaystyle{\lfloor\frac{w_{\alpha}+\rho}{m}\rfloor((\frac{\alpha+2}{2})\nu-m)-\lfloor\frac{w_{\alpha-1}+\rho}{m}\rfloor(\frac{\alpha\nu}{2}).}
\end{array}
$

\noindent By \eqref{4.1 2}, we have
$
\displaystyle{\sum_{j={\alpha+1}}^{m-1}(\lfloor\frac{w_j+\rho}{m}\rfloor-\lfloor\frac{w_{j-1}+\rho}{m}\rfloor)(j\nu-m)}\geq \displaystyle{(\lfloor\frac{w_{m-1}+\rho}{m}\rfloor-\lfloor\frac{w_{\alpha}+\rho}{m}\rfloor)((\alpha+1)\nu-m).}
$

\noindent Since $w_{m-1}\geq w_{\alpha-1}+w_{\alpha}$, by Lemma \ref{x}, it follows that $\lfloor\frac{w_{m-1}+\rho}{m}\rfloor\geq\lfloor\frac{w_{\alpha-1}+\rho}{m}\rfloor+\lfloor\frac{w_{\alpha}+\rho}{m}\rfloor-1.$ Let $x=\lfloor\frac{w_{m-1}+\rho}{m}\rfloor\-\lfloor\frac{w_{\alpha-1}+\rho}{m}\rfloor-\lfloor\frac{w_{\alpha}+\rho}{m}\rfloor$. Then, $\lfloor\frac{w_{\alpha-1}+\rho}{m}\rfloor+\lfloor\frac{w_{\alpha}+\rho}{m}\rfloor=q-x$ and $x\geq -1$. Now using $\rho\geq \frac{(3-q)\alpha m}{2\alpha+6}$ and $(\frac{\alpha+3}{3})\nu\geq m$, we get 
\begin{equation*}
\begin{array}{lllll}
\displaystyle{\sum_{j=1}^{m-1}(\lfloor\frac{w_j+\rho}{m}\rfloor-\lfloor\frac{w_{j-1}+\rho}{m}\rfloor)(j\nu-m)+\rho
}
&\geq \displaystyle{\lfloor\frac{w_{\alpha}+\rho}{m}\rfloor((\frac{\alpha+2}{2})\nu-m)-\lfloor\frac{w_{\alpha-1}+\rho}{m}\rfloor(\frac{\alpha\nu}{2})}\\&+(\displaystyle{\lfloor\frac{w_{m-1}+\rho}{m}\rfloor-\lfloor\frac{w_{\alpha}+\rho}{m}\rfloor)((\alpha+1)\nu-m)+\rho
}\\ \\
&=\displaystyle{\lfloor\frac{w_{\alpha-1}+\rho}{m}\rfloor(\frac{-\alpha\nu}{2}+(\alpha+1)\nu-m-((\alpha+1)\nu-m))}\\&+\displaystyle{\lfloor\frac{w_{\alpha}+\rho}{m}\rfloor((\frac{\alpha+2}{2})\nu-m)}+(\displaystyle{\lfloor\frac{w_{m-1}+\rho}{m}\rfloor-\lfloor\frac{w_{\alpha}+\rho}{m}\rfloor)((\alpha+1)\nu-m)}\\&+\rho\\ \\

&= (\displaystyle{\lfloor\frac{w_{\alpha-1}+\rho}{m}\rfloor+\lfloor\frac{w_{\alpha}+\rho}{m}\rfloor)((\frac{\alpha+2}{2})\nu-m)}\\ &+\displaystyle{(\lfloor\frac{w_{m-1}+\rho}{m}\rfloor-\lfloor\frac{w_{\alpha-1}+\rho}{m}\rfloor-\lfloor\frac{w_{\alpha}+\rho}{m}\rfloor)((\alpha+1)\nu-m)+\rho}\\ \\

&= (q-x)((\frac{\alpha+2}{2})\nu-m)+x((\alpha+1)\nu-m)+\rho\\ \\

&= \nu(q+\frac{q\alpha}{2}+\frac{x\alpha}{2})-qm+\rho\\ \\
&\geq \nu(q+\frac{q\alpha}{2}-\frac{\alpha}{2})-qm+ \frac{(3-q)\alpha m}{2\alpha+6}\\ \\

&=\nu(q+\frac{q\alpha}{2}-\frac{\alpha}{2})-m(\frac{q(2\alpha+6)+(q-3)\alpha}{2\alpha+6})\\ \\

&=\nu(q+\frac{q\alpha}{2}-\frac{\alpha}{2})-m(\frac{3q}{\alpha+3}+\frac{3q\alpha}{2(\alpha+3)}-\frac{3\alpha}{2(\alpha+3)})\\ \\

&= (q+\frac{q\alpha}{2}-\frac{\alpha}{2})(\frac{3}{\alpha+3})((\frac{\alpha+3}{3})\nu-m)\\ \\

&\geq 0.
\end{array}
\end{equation*}

\noindent Using  Proposition \ref{prop2}, we get that $S$ satisfies Wilf's Conjecture. $  \quad \square$

\end{demostracion}

\begin{exemple}\label{ex2} {\rm Consider the following numerical semigroup $S=<22,23,25,27,29,31,33>$. Note that $3\nu<m$. We have $w_6=33
,$ $w_7=46$ and $w_{m-1}=87$ and i.e. $w_{m-1}\geq w_6+w_7 $. Moreover, $ (\frac{\alpha+3}{3})\nu=(\frac{7+3}{3})7\geq 22=m$, thus the conditions of Theorem \ref{thm3} are valid. Hence, $S$ satisfies Wilf's Conjecture.}$\quad \square$
\end{exemple}
\medskip
\noindent Th following Corollary \ref{cory} is an extension for Corollary \ref{cor1} using Theorems \ref{thm2} and \ref{thm3}.
\begin{corolario} \label{cory}{\rm Let $S$ be a numerical semigroup with multiplicity $m$ and embedding dimension $\nu$. Suppose that $m-\nu\geq \frac{\alpha(\alpha-1)}{2}-1$ for some $7\leq\alpha\leq m-2$. If $(2+\frac{\alpha-3}{q})\nu\geq m$, then $S$ satisfies Wilf's Conjecture.}
\end{corolario}

\begin{demostracion}{.} Since $m-\nu\geq \frac{\alpha(\alpha-1)}{2}-1$, then by Lemma \ref{lem2} we have $w_{m-1}\geq w_1+w_{\alpha}$ or $w_{m-1}\geq w_{\alpha-2}+w_{\alpha-1}$. Suppose that $w_{m-1}\geq w_1+w_{\alpha}$. Since $(2+\frac{\alpha-3}{q})\nu\geq m$, by applying Theorem \ref{thm2}, $S$ satisfies wilf's Conjecture. Now suppose that $w_{m-1}\geq w_{\alpha-2}+w_{\alpha-1}$. We may assume that $q\geq 4$ ($q\leq3$ is solved \cite{Kaplan1}, \cite{Eliahou}). Then for $\alpha\geq 7$ we have, ($\frac{\alpha-1+3}{3})\nu\geq(2+\frac{\alpha-3}{q})\nu$. Consequently,  ($\frac{\alpha-1+3}{3})\nu\geq m$. Now by applying Theorem \ref{thm3}, $S$ satisfies Wilf's Conjecture. $\quad \square$
\end{demostracion}

\medskip
\noindent As a direct consequence of Theorem \ref{thm3}, we get the following Corollary.
\begin{corolario}\label{cor4} {\rm Let $S$ be a numerical semigroup with multiplicity $m$ and embedding dimension $\nu$. Let $w_0=0<w_1<w_2<\ldots<w_{m-1}$ be the elements of Ap($S,m)$. 
Suppose that $ w_{m-1}\geq w_{\alpha-1}+w_{\alpha}$ for some $1<\alpha<m-1$. If $m\leq \frac{4(\alpha+3)}{3}$, then $S$ satisfies Wilf's Conjecture.}
\end{corolario}
\begin{demostracion}{.}
If $\nu<4$, then $S$ satisfies Wilf's Conjecture (\cite{Dobbs1}). Hence, we can suppose that $\nu\geq 4$. Thus, $(\frac{\alpha+3}{3})(\nu)\geq \frac{4(\alpha+3)}{3}\geq m$. By applying Theorem \ref{thm3} $S$ satisfies Wilf's Conjecture. $\quad \square$
\end{demostracion}

\end{document}